\documentclass[11pt,a4paper]{article}
\usepackage{amsthm,mathrsfs,amsmath,amssymb}
\usepackage{dsfont}
\usepackage{comment}
\usepackage{indentfirst}
\usepackage{cite}
\usepackage[colorlinks=true]{hyperref}
\hypersetup{urlcolor=blue, citecolor=red}
\numberwithin{equation}{section}
\allowdisplaybreaks
\makeatother
\date{}

\newtheorem{thm}{Theorem}[section]
\newtheorem{lem}[thm]{Lemma}
\newtheorem{prop}[thm]{Proposition}

\newtheorem{xrem}[thm]{\bf Remark}
\newenvironment{prf}{\noindent {\bf Proof.}\rm}{\qed}
\newenvironment{prf1}{\noindent {\bf Proof of Theorem \ref{T1}.}\rm}{\qed}
\newenvironment{prf2}{\noindent {\bf Proof of Theorem \ref{T2}.}\rm}{\qed}
\makeatother

\title{\bf Normalized solution for lower critical nonhomogeneous Choquard equation\footnote{Supported by National Natural Science Foundation of China (No.~12301143, No.~12261079), Natural Science Foundation of Gansu Province, China (No.~24JRRA141) and NWNU-LKQN2022-02.}}

\author{{Yong-Yong Li\footnote{Corresponding author. E-mail address: mathliyy518@163.com (Y.-Y. Li).}}\\
{\small\emph{College of Mathematics and Statistics, Northwest Normal University, Lanzhou, {\rm730070}, China,}}\\}

\begin{document}
\maketitle

\baselineskip 16pt

\setcounter{page}{1}


\begin{quote}
  {\bf Abstract:}\ In this paper, we study the lower critical nonhomogeneous Choquard equation
  \begin{align*}
  \begin{cases}
  \displaystyle -\Delta u+\lambda u=\left(I_\alpha\ast|u|^{2_*^\alpha}\right)|u|^{2_*^\alpha-2}u+h(x) &\mbox{in}\ \mathbb{R}^N, \\
  \displaystyle \int_{\mathbb{R}^N}|u|^2dx=a^2,
  \end{cases}
  \end{align*}
  where $N\geq3$, $a>0$, $\lambda\in\mathbb{R}$ appears as an unknown Lagrange multiplier, $2_*^\alpha=\frac {N+\alpha}N$, $I_\alpha$ is the Riesz potential of order $\alpha\in(0,N)$ and $h\in L^2(\mathbb{R}^N)$ satisfies some mild conditions.
  If $h(\cdot)>0$, based on a rather careful analysis for the minimizing sequence, we get a positive ground state normalized solution by applying a constrained minimization argument on $L^2$ sphere. If $h(\cdot)\leq0$, we verify the nonexistence of nonnegative minimizer for the~least~energy.

\end{quote}


\section{Introduction and main results}

In recent decades, many scholars studied various variational problems for the Choquard equation
\begin{align}\label{CH1}
-\Delta u+\lambda u=\left(I_\alpha\ast|u|^p\right)|u|^{p-2}u\ \ \ &\mbox{in}\ \mathbb{R}^N,
\end{align}
where $N\geq3$, $2_*^\alpha:=\frac{N+\alpha}N\leq p\leq 2_\alpha^*:=\frac{N+\alpha}{N-2}$,
$I_\alpha$ is the Riesz potential of order $\alpha\in(0,N)$ defined by
$$I_\alpha(\cdot)=\frac{A_\alpha}{|\cdot|^{N-\alpha}}=\frac{\Gamma(\frac{N-\alpha}2)}{{2^\alpha\pi^{N/2}\Gamma(\frac\alpha2)|\cdot|^{N-\alpha}}},$$
$\Gamma$ is Gamma function and $*$ is convolution in $\mathbb{R}^N$.
The endpoints $2_\alpha^*$ and $2_*^\alpha$ are usually called as upper and lower critical exponents due to the following Hardy-Littlewood-Sobolev inequality, respectively.
\begin{prop}\label{P1}
(Theorem 4.3, \cite{LL}). Let $r,s>1$ and $0<\mu<N$ satisfy $\frac1r+\frac1s+\frac\mu N=2$.
Then there exists a sharp constant $C(N,\mu,r)>0$ such that, for all $u\in L^r(\mathbb{R}^N)$ and $v\in L^s(\mathbb{R}^N)$, there holds
\begin{align}\label{HLS}
\left|\int_{\mathbb{\mathbb{R}}^N}\int_{\mathbb{\mathbb{R}}^N}
\frac{u(x)v(y)}{|x-y|^{\mu}}dxdy\right|
\leq C(N,\mu,r)|u|_r|v|_s.
\end{align}
\end{prop}

\vspace{0.10cm}
This inequality is sometimes referred to as a weak version of the following Young inequality:

\begin{prop}\label{P2}
(Theorem 4.2, \cite{LL}). Let $r, s, t\geq1$ satisfy $\frac1r+\frac1s+\frac1t=2$.
Then there exists a sharp constant $C(N,r,s)>0$ such that, for all $u\in L^r(\mathbb{R}^N)$, $w\in L^s(\mathbb{R}^N)$ and $v\in L^t(\mathbb{R}^N)$, there holds
\begin{align}\label{Y}
\left|\int_{\mathbb{R}^N}\int_{\mathbb{R}^N}u(x)w(x-y)v(y)dxdy\right|
\leq C(N,r,s)|u|_r|w|_s|v|_t.
\end{align}
\end{prop}

As we know, Eq.~\eqref{CH1} has a physical prototype (i.e.~the Hartree type evolution equation), which firstly appeared as a model of static polaron \cite{P},
then Choquard introduced it as a model of electron well \cite{L} and Penrose proposed it as a gravitational model in quantum gravity \cite{P1}.
In mathematical contexts, fixed frequency problems (i.e.~$\lambda\in\mathbb{R}$ is given) for Eq.~\eqref{CH1} were well studied,
see \cite{CVZ,LLT,MV,MV1,VX} and the references therein.
In particular, since $2_*^\alpha$ provokes a new kind of loss of compactness due to the invariance of scaling action and the authors in \cite{MV1} proved that Eq.~\eqref{CH1} has nontrivial solution iff $p\in(2_*^\alpha,2_\alpha^*)$, the lower critical Choquard equations with noncritical perturbations have been attracting  widespread attention, see e.g. \cite{CVZ,LLT,VX}.
While, for Eq.~\eqref{CH1}, the prescribed mass problems
\begin{align}\label{CH2}
\begin{cases}
\displaystyle -\Delta u+\lambda u=\gamma\left(I_\alpha\ast|u|^p\right)|u|^{p-2}u+g(x,u)  &\mbox{in}\ \mathbb{R}^N, \\
\displaystyle |u|_2=a
\end{cases}
\end{align}
are seldom studied, where $\lambda\in\mathbb{R}$ appears as Lagrange multiplier, $g\in C(\mathbb{R}^N\times\mathbb{R},\mathbb{R})$ and $\gamma, a>0$, which is closely related to the Bose-Einstein condensation, such solutions are called as normalized~solutions.

For the case $p=2_*^\alpha$ of Eq.~\eqref{CH2}, there is no normalized solution if $g=0$.
Yao et al. \cite{YCRS} studied the case that $p=2_*^\alpha$, $g(x,u)=\mu|u|^{q-2}u$, $q\in\big(2,\frac{2N}{N-2}\big]$ and $\mu>0$ of Eq.~\eqref{CH2},
under some assumptions on $\gamma, \mu, a$, they proved ground state normalized solution when $q\in(2,\bar{q})\cup(\bar{q},\frac{2N}{N-2}]$, also showed the nonexistence of normalized solution when $q=\bar{q}$, where $\bar{q}=2+\frac4N$ is the so-called L$^2$-critical exponent. Afterwards, Li et al. \cite{LBT} improved some results in \cite{YCRS}. 
For more studies on existence and multiplicity of normalized solutions for the Choquard equations, we refer to \cite{BLL,LSZ,YHSW,Y,YSY,YCT} and the references therein.
Having recalled the above results, we find a question whether Eq.~\eqref{CH2} with sublinear perturbation has normalized solution left.
Motivated by \cite{CZ}, where Chen and Zou studied the Schr\"odinger equation
\begin{align*}
\begin{cases}
\displaystyle -\Delta u+\lambda u=|u|^{q-2}u+h(x)\ \ \ \mbox{in}\ \mathbb{R}^N,\\
\displaystyle |u|_2=a,
\end{cases}
\end{align*}
when $q\in(2,\bar{q})\cup(\bar{q},\frac{2N}{N-2})$ and $h$ satisfies some mild hypotheses they showed the existence of normalized solution,
we study the nonhomogeneous Choquard equation with lower critical exponent here, namely,
\begin{align}\label{CH3}
\begin{cases}
\displaystyle -\Delta u+\lambda u=\left(I_\alpha\ast|u|^{2_*^\alpha}\right)|u|^{2_*^\alpha-2}u+h(x) &\mbox{in}\ \mathbb{R}^N, \\
\displaystyle |u|_2=a,
\end{cases}
\end{align}
where $N\geq3$, $\alpha\in(0,N)$, $a>0$, $\lambda\in\mathbb{R}$ is unknown Lagrange multiplier and $h$ may satisfy the~condition
\begin{itemize}
  \item [($\mathcal{H}$)] there exist
  $\beta\in\left(\frac N2, \min\left\{\frac N2+2,N\right\}\right)$ and $\bar{r}>0$ such that $h(x)\geq|x|^{-\beta}$ for all $|x|\geq \bar{r}$.
\end{itemize}
Technically, normalized solutions of Eq.~\eqref{CH3} can be catched as critical points of the energy functional
\begin{align*}
\mathcal{J}(u)=\frac12\int_{\mathbb{R}^N}|\nabla u|^2dx
-\frac1{2\cdot2_*^\alpha}\int_{\mathbb{R}^N}\left(I_\alpha\ast|u|^{2_*^\alpha}\right)|u|^{2_*^\alpha}dx-\int_{\mathbb{R}^N}h(x)udx
\end{align*}
on the constraint
\begin{align*}
\mathcal{S}(a):=\left\{u\in H^1(\mathbb{R}^N):|u|_2=a\right\}.
\end{align*}
For any $a>0$, we define that $u$ is a ground state normalized solution of Eq.~\eqref{CH3} if $u\in\mathcal{S}(a)$ satisfies
\begin{align*}
\mathcal{J}|_{\mathcal{S}(a)}'(u)=0\ \ \ \ \mbox{and}\ \ \ \ \mathcal{J}(u)=\inf\left\{\mathcal{J}(v):v\in \mathcal{S}(a), \ \mathcal{J}|_{\mathcal{S}(a)}'(v)=0\right\}.
\end{align*}
To search for ground state normalized solution of Eq.~\eqref{CH3}, we consider the minimization problem
\begin{align*}
m(a):=\inf_{u\in\mathcal{S}(a)}\mathcal{J}(u).
\end{align*}

Now we are prepared to state our main result:
\begin{thm}\label{T1}
Assume $0< h(\cdot)\in L^2(\mathbb{R}^N)$ and $(\mathcal{H})$ holds, then Eq.~\eqref{CH3} admits a positive ground state normalized solution, with the corresponding Lagrange multiplier $\lambda_a>\frac{a^{2(2_*^\alpha-1)}}{2_*^\alpha S_\alpha^{2_*^\alpha}}-\frac{|h|_2}{a}$.
\end{thm}

\begin{thm}\label{T2}
If $h(\cdot)\leq 0$ and $h\in L^2(\mathbb{R}^N)$, then $m(a)$ has no nonnegative minimizer for any $a>0$.
\end{thm}

\begin{xrem}\label{R1.2}
In \cite{Y}, Ye investigated Eq.~\eqref{CH3} with $h=0$ and showed that the minimization problem
\begin{align}\label{R0}
m_\infty(a):=\inf_{u\in\mathcal{S}(a)}\left[\frac12\int_{\mathbb{R}^N}|\nabla u|^2dx
-\frac1{2\cdot2_*^\alpha}\int_{\mathbb{R}^N}\left(I_\alpha\ast|u|^{2_*^\alpha}\right)|u|^{2_*^\alpha}dx\right]
=-\frac1{2\cdot2_*^\alpha}\left(\frac{a^2}{S_\alpha}\right)^{2_*^\alpha}
\end{align}
has no minimizer. Theorem \ref{T1} is a first counterpart of the result \cite{CZ} on Sch\"{o}dinger equations for the Choquard equations.
However, comparing to \cite{CZ}, the limit problem \eqref{R0} of Eq.~\eqref{CH3} has no minimizer and the Lions lemma is invalid for the lower critical term,
thus the methods there cannot be directly used to prove our result,
we need create new tricks and impose $(\mathcal{H})$ to restore the lack of compactness.
Moreover, we mention that the subadditivity inequality of least energy as in \cite{LBT,YCRS} may not hold here.
This study complements the nonexistence result of lower critical Choquard equation in \cite{Y}, meanwhile
extends the works \cite{YCRS,LBT} from the case of superlinear perturbation to the case of sublinear perturbation.
\end{xrem}

\noindent {\bf Notations:}
\begin{enumerate}
  \setlength{\itemsep}{-0.5pt}
  \item [$\spadesuit$] $L^p(\mathbb{R}^N)$ is the usual Lebesgue space with the norm $|u|_p=\left(\int_{\mathbb{R}^N}|u|^pdx\right)^{\frac1p}$ for $p\in[1,+\infty)$.
  \item [$\spadesuit$] $C_0^\infty(\mathbb{R}^N)$ consists of infinitely times differentiable functions with compact support in $\mathbb{R}^N$.
  \item [$\spadesuit$] $H^1(\mathbb{R}^N)=\left\{u\in L^2(\mathbb{R}^N):|\nabla u|\in L^2(\mathbb{R}^N)\right\}$ endowed with the norm
  $\|u\|=(\int_{\mathbb{R}^N}|\nabla u|^2+u^2dx)^{\frac12}$.
  \item [$\spadesuit$] $\mathbb{B}_r(y)=\left\{x\in\mathbb{R}^N:|x-y|<r\right\}$ and $\mathbb{B}_r^c(y)=\mathbb{R}^N\backslash\mathbb{B}_r(y)$  for $r>0$ and $y\in\mathbb{R}^N$.
  \item [$\spadesuit$] The best constant $S_\alpha:=\inf\left\{|u|_2^2:\int_{\mathbb{R}^N}(I_\alpha\ast|u|^{2_*^\alpha})|u|^{2_*^\alpha}dx=1\right\}$.
  \item [$\spadesuit$] $o(1)$ denotes a quantity tending to $0$ as $n\rightarrow\infty$.
\end{enumerate}

\section{Proof of Theorems \ref{T1} and \ref{T2}}

Before completing the proof of Theorem \ref{T1}, we establish some preliminaries for the least energy.

\begin{lem}\label{L1}
(1) The functional $\mathcal{J}$ is coercive on $\mathcal{S}(a)$ and $-\infty<m(a)< m_\infty(a)$ for any $a>0$;
\begin{itemize}
  \item [(2)] $m(a)$ is continuous with respect to $a$ on $(0,+\infty)$;
  \item [(3)] $m\big(\sqrt{a_1^2+a_2^2}\big)\leq m(a_1)+m_\infty(a_2)$ for any $a_1, a_2>0$.
\end{itemize}
\end{lem}

\begin{prf}
(1) By the definition of $S_\alpha$ and the H$\ddot{\mbox{o}}$lder inequality, we derive that $\mathcal{J}$ is coercive on $\mathcal{S}(a)$.
Naturally, $m(a)>-\infty$. Recall that $S_\alpha$ is achieved by $u_0(x)= A{(1+|x|^2)^{-\frac N2}}$, where $A>0$ is fixed such that
$$
\int_{\mathbb{R}^N}\left(I_\alpha\ast|u_0|^{2_*^\alpha}\right)|u_0|^{2_*^\alpha}dx=1.
$$
Set $u_a^\varepsilon(\cdot)=\frac{a\varepsilon^{\frac N2}u_0(\varepsilon\cdot)}{|u_0|_2}$ for $\varepsilon>0$.
Clearly, $u_a^\varepsilon\in\mathcal{S}(a)$.
Then it follows from $(\mathcal{H})$ that, once $\varepsilon\in(0,1)$,
\begin{align}\label{R4}
m(a)\leq\mathcal{J}(u_a^\varepsilon)
&\leq\frac{a^2\varepsilon^2}{2|u_0|_2^2}\int_{\mathbb{R}^N}|\nabla u_0|^2dx
-\frac1{2\cdot2_*^\alpha}\left(\frac a{|u_0|_2}\right)^{2\cdot2_*^\alpha}
-\frac{a\varepsilon^{\frac N2}}{|u_0|_2}\int_{\mathbb{B}_{\frac {\bar{r}}\varepsilon}^c(0)}\frac{u_0(\varepsilon x)}{|x|^\beta}dx \nonumber\\
&\leq \frac{a^2\varepsilon^2}{2|u_0|_2^2}\int_{\mathbb{R}^N}|\nabla u_0|^2dx
-\frac{a\varepsilon^{\beta-\frac N2}}{|u_0|_2}\int_{\mathbb{B}_{\bar{r}}^c(0)}\frac{u_0(x)}{|x|^\beta}dx
-\frac1{2\cdot2_*^\alpha}\left(\frac {a^2}{S_\alpha}\right)^{2_*^\alpha},
\end{align}
by which, noting $0<\beta-\frac N2<2$, for $\varepsilon>0$ small enough we deduce $m(a)<m_\infty(a)$. Hence, (1) is true.

(2) For any $a>0$ and sequence $\{a_n\}\subset (0,+\infty)$ satisfying $a_n\rightarrow a$, by the definition of $m(a_n)$, for any $\varepsilon>0$, there exists $u_n\in \mathcal{S}(a_n)$ such that $\mathcal{J}(u_n)<m(a_n)+\varepsilon$ for all $n$. Clearly, $\{u_n\}$ is bounded in $H^1(\mathbb{R}^N)$. Set $\tilde{u}_n=\frac{au_n}{|u_n|_2}$, there holds
\begin{align*}
m(a)\leq\mathcal{J}(\tilde{u}_n)=\mathcal{J}(u_n)+o(1)\leq m(a_n)+\varepsilon+o(1).
\end{align*}
Similarly, by the definition of $m(a)$, for any positive sequence $\{\varepsilon_n\}$ with $\varepsilon_n\rightarrow0$, there exists a sequence $\{v_n\}\subset\mathcal{S}(a)$ such that $\mathcal{J}(v_n)<m(a)+\varepsilon_n$ for all $n$. Set $\tilde{v}_n=\frac{v_na_n}{|v_n|_2}$.
Since $\{\|v_n\|\}$ is bounded, we~deduce
\begin{align*}
m(a_n)\leq \mathcal{J}(\tilde{v}_n)=\mathcal{J}(v_n)+o(1)\leq m(a)+\varepsilon_n+o(1).
\end{align*}
By combining the above two inequalities, we may conclude $m(a_n)\rightarrow m(a)$. That is, (2) is certified.

(3) For any $\varepsilon>0$, by $\overline{C_0^\infty(\mathbb{R}^N)}= H^1(\mathbb{R}^N)$,
there exist some $\varphi_\varepsilon, \psi_\varepsilon \in C_0^\infty(\mathbb{R}^N)$ such that
\begin{align*}
\begin{cases}
\vspace{0.10cm}
\displaystyle \varphi_\varepsilon\in\mathcal{S}(a_1)\ \mbox{and} \ \mathcal{J}(\varphi_\varepsilon)<m(a_1)+\frac\varepsilon2,\\
\displaystyle \psi_\varepsilon\in\mathcal{S}(a_2)\ \mbox{and} \ \mathcal{J}_\infty(\psi_\varepsilon)<m_\infty(a_2)+\frac\varepsilon2.
\end{cases}
\end{align*}
Define $\omega_{\varepsilon,n}(\cdot)=\varphi_\varepsilon(\cdot)-\psi_\varepsilon(\cdot-n\mathbf{e}_1)$, where $\mathbf{e}_1=(1,0,\cdots,0)\in\mathbb{R}^N$.
For large $n$, by $\mbox{spt}~\varphi_\varepsilon\cap \mbox{spt}~\psi_\varepsilon(\cdot-n\mathbf{e}_1)=\emptyset$, we deduce $\omega_{\varepsilon,n}\in \mathcal{S}(\sqrt{a_1^2+a_2^2})$ and then
\begin{align*}
m\Big(\sqrt{a_1^2+a_2^2}\Big)\leq\mathcal{J}(\omega_{\varepsilon,n})=\mathcal{J}(\varphi_\varepsilon)+\mathcal{J}(\psi_\varepsilon(\cdot-n\mathbf{e}_1))
-\frac1{2_*^\alpha}\int_{\mathbb{R}^N}\left(I_\alpha\ast|\varphi_\varepsilon|^{2_*^\alpha}\right)\left|\psi_\varepsilon(\cdot-n\mathbf{e}_1)\right|^{2_*^\alpha}dx.
\end{align*}
Since $\psi_\varepsilon(\cdot-n\mathbf{e}_1)\xrightarrow{n}0$ a.e.~in $\mathbb{R}^N$,
then $\psi_\varepsilon(\cdot-n\mathbf{e}_1)\overset{n}\rightharpoonup0$ in $L^2(\mathbb{R}^N)$
and $|\psi_\varepsilon(\cdot-n\mathbf{e}_1)|^{2_*^\alpha}\overset{n}\rightharpoonup0$ in $L^{\frac{2N}{N+\alpha}}(\mathbb{R}^N)$ up to a subsequence.
Further, it follows from $h\in L^2(\mathbb{R}^N)$ and $I_\alpha\ast|\varphi_\varepsilon|^{2_*^\alpha}\in L^{\frac{2N}{N-\alpha}}(\mathbb{R}^N)$ that as $n\rightarrow\infty$,
\begin{align*}
\int_{\mathbb{R}^N}h(x)\psi_\varepsilon(\cdot-n\mathbf{e}_1)dx\rightarrow0
 \end{align*}
and
\begin{align*}
\int_{\mathbb{R}^N}\left(I_\alpha\ast|\varphi_\varepsilon|^{2_*^\alpha}\right)|\psi_\varepsilon(\cdot-n\mathbf{e}_1)|^{2_*^\alpha}dx\rightarrow0.
\end{align*}
Consequently, there results
\begin{align*}
m\Big(\sqrt{a_1^2+a_2^2}\Big)
&\leq\limsup_{n\to\infty}\big[\mathcal{J}(\varphi_\varepsilon)+\mathcal{J}(\psi_\varepsilon(\cdot-n\mathbf{e}_1))\big]\\
&=\mathcal{J}(\varphi_\varepsilon)+\mathcal{J}_\infty(\psi_\varepsilon)\\
&< m(a_1)+m_\infty(a_2)+\varepsilon,
\end{align*}
which together with the arbitrariness of $\varepsilon>0$ implies that (3) holds.
Thus this lemma is proved.
\end{prf}\\

\begin{prf1}
For any $a>0$, take a minimizing sequence $\{u_n\}\subset\mathcal{S}(a)$ for $m(a)$.
By Lemma \ref{L1}-(1), $\{u_n\}$ is bounded in $H^1(\mathbb{R}^N)$.
Thanks to the the concentration-compactness principle (see \cite[Lemma III.1]{L1}),
in the sense of subsequence, there product the following three possibilities:
\begin{enumerate}
  \item [(a)] vanishing:
  $\displaystyle \lim_{n\to\infty}\sup_{y\in\mathbb{R}^N}\int_{\mathbb{B}_r(y)}|u_n(x)|^2dx=0$ for any $r>0$;
  \item [(b)] dichotomy: there exist $a_1\in (0,a)$ and sequences $\{u_{n,1}\}, \{u_{n,2}\}$ bounded in $H^1(\mathbb{R}^N)$ such that
  \begin{align}\label{R8}
  \begin{cases}
  \displaystyle \left|u_n-u_{n,1}-u_{n,2}\right|_2\rightarrow 0, \\
  \displaystyle \left|u_{n,1}\right|_2^2\rightarrow a_1^2\ \mbox{and}\ \left|u_{n,2}\right|_2^2\rightarrow a^2-a_1^2, \\
  \displaystyle \mbox{dist}\left(\mbox{spt}~u_{n,1},\mbox{spt}~u_{n,2}\right)\rightarrow+\infty, \\
  \displaystyle \liminf\limits_{n\rightarrow\infty}\int_{\mathbb{R}^N}|\nabla u_n|^2-|\nabla u_{n,1}|^2-|\nabla u_{n,2}|^2dx\geq0,
  \end{cases}
  \end{align}
  we point out $u_{n,1}(\cdot)=\varphi\big(\frac{\cdot-y_n}{R_0}\big)u_n(\cdot)$ for some $R_0>0$ and $u_{n,2}(\cdot)=\psi\big(\frac{\cdot-y_n}{R_n}\big)u_n(\cdot)$ with $R_n\rightarrow+\infty$, where $\{y_n\}\subset\mathbb{R}^N$, the cut-off functions $\varphi, \psi\in C^\infty(\mathbb{R}^N,[0,1])$ such that $\varphi(x)=1$ if $|x|\leq1$, $\varphi(x)=0$ if $|x|\geq2$ and $\psi(x)=0$ if $|x|\leq1$, $\psi(x)=1$ if $|x|\geq2$;
  \item [(c)] compactness: there exists $\{z_n\}\subset\mathbb{R}^N$ such that, for any $\varepsilon>0$, there exists some $R>0$ satisfying
  \begin{align}\label{R1}
  \int_{\mathbb{B}_R(z_n)}|u_n(x)|^2dx\geq a^2-\varepsilon.
  \end{align}
\end{enumerate}

If vanishing occurs, then $u_n\to0$ in $L^p(\mathbb{R}^N)$ for any $p\in(2,\frac{2N}{N-2})$, therefore $u_n\rightharpoonup0$ in $L^2(\mathbb{R}^N)$ up to a subsequence. Further, due to Lemma \ref{L1}-(1) and $h\in L^2(\mathbb{R}^N)$, there holds
\begin{align*}
m_\infty(a)
>m(a)
&=\lim_{n\to\infty}\left[\frac12\int_{\mathbb{R}^N}|\nabla u_n|^2dx
-\frac1{2\cdot2_*^\alpha}\int_{\mathbb{R}^N}\left(I_\alpha\ast|u_n|^{2_*^\alpha}\right)|u_n|^{2_*^\alpha}dx\right]\\
&\geq -\frac1{2\cdot2_*^\alpha}\left(\frac {a^2}{S_\alpha}\right)^{2_*^\alpha},
\end{align*}
a contradiction. We claim dichotomy cannot occur.
Indeed, if $\{y_n\}$ is unbounded, assume $|y_n|\rightarrow+\infty$, then $u_{n,1}(x)\rightarrow0$ a.e.~in $\mathbb{R}^N$.
It is easy to show
\begin{align}\label{R11}
\int_{\mathbb{R}^N}h(x)u_{n,1}dx\to0.
\end{align}
Moreover, due to \cite[Lemma 3.3]{CVZ}, the functional
\begin{align*}
\|\cdot\|_*:=\left[\int_{\mathbb{R}^N}\big(I_\alpha\ast|\cdot|^{2_*^\alpha}\big)|\cdot|^{2_*^\alpha}dx\right]^{\frac1{2\cdot2_*^\alpha}}
\end{align*}
is a norm in $L^2(\mathbb{R}^N)$. Then there holds
\begin{align}\label{R9}
&\left|\|u_n\|_*-\|u_{n,1}+u_{n,2}\|_*\right| \nonumber\\
&\quad\leq\|u_n-u_{n,1}-u_{n,2}\|_* \nonumber\\
&\quad\leq S_\alpha^{-\frac12}\left|u_n-u_{n,1}-u_{n,2}\right|_2\xrightarrow{n}0.
\end{align}
Set $f(\cdot)=|\cdot|^{\alpha-N}$ and $f_\delta(\cdot)=f(\cdot)\chi_{\{f\geq\delta\}}(\cdot)$ for $\delta>0$, where $\chi$ denotes characteristic function.
Since $\mbox{spt}~u_{n,1}\subset\mathbb{B}_{2R_0}(y_n)$, due to \eqref{HLS} and the H$\ddot{\mbox{o}}$lder inequality,
there exists some $C(N,\alpha)>0$ such that
\begin{align*}
\int_{\mathbb{R}^N}&\int_{\mathbb{R}^N}\frac{|u_{n,1}(x)|^{2_*^\alpha}|u_{n,2}(y)|^{2_*^\alpha}\chi_{\{f<\delta\}}(x-y)}{|x-y|^{N-\alpha}}dxdy \\
&\leq\delta^{\frac14}\int_{\mathbb{R}^N}\int_{\mathbb{R}^N}\frac{|u_{n,1}(x)|^{2_*^\alpha}|u_{n,2}(y)|^{2_*^\alpha}}{|x-y|^{\frac{3(N-\alpha)}4}}dxdy \\
&\leq \delta^{\frac14}C(N,\alpha)|u_{n,2}|_2^{2_*^\alpha}\left(\int_{\mathbb{R}^N}|u_{n,1}|^{\frac{4(N+\alpha)}{3N+\alpha}}dx\right)^{\frac{3N+\alpha}{4N}} \\
&\leq \delta^{\frac14} C(N,\alpha) a^{2\cdot2_*^\alpha} \left|\mbox{spt}~u_{n,1}\right|^{\frac{N-\alpha}{4N}}.
\end{align*}
Denote $d_n=\mbox{dist}\left(\mbox{spt}~u_{n,1},\mbox{spt}~u_{n,2}\right)$ and
\begin{align*}
\sigma=
\begin{cases}
\displaystyle 1  &\mbox{if} \ 0<\alpha<\frac{2N}{N-2}, \\
\displaystyle \frac{N^2}{2N^2-(N+\alpha)(N-2)} &\mbox{if} \ \frac{2N}{N-2}\leq\alpha<N.
\end{cases}
\end{align*}
Obviously, $\sigma\geq1$ and $\frac{2\cdot2_*^\alpha\sigma}{2\sigma-1}\in\big(2,\frac{2N}{N-2}\big]$ for any $\alpha\in(0,N)$.
It is easy to show $f_\delta\in L^\sigma(\mathbb{R}^N)$ and then $\left|f_\delta\chi_{\{|x|\geq d_n\}}\right|_\sigma\rightarrow0$ since $d_n\rightarrow+\infty$.
By Young's inequality \eqref{Y} and the Sobolev inequality, there exists some $\bar{C}(N,\alpha)>0$ such that
\begin{align*}
\int_{\mathbb{R}^N}&\int_{\mathbb{R}^N}{|u_{n,1}(x)|^{2_*^\alpha}|u_{n,2}(y)|^{2_*^\alpha}}{f_\delta(x-y)}dxdy \\
&=\int_{\mathbb{R}^N}\int_{\mathbb{R}^N}{|u_{n,1}(x)|^{2_*^\alpha}|u_{n,2}(y)|^{2_*^\alpha}}{f_\delta(x-y)\chi_{\{|x-y|\geq d_n\}}(x-y)}dxdy \\
&\leq \bar{C}(N,\alpha)|u_{n,1}|_{\frac{2\cdot2_*^\alpha\sigma}{2\sigma-1}}^{2_*^\alpha}
|u_{n,2}|_{\frac{2\cdot2_*^\alpha\sigma}{2\sigma-1}}^{2_*^\alpha}\left|f_\delta\chi_{\{|x|\geq d_n\}}\right|_\sigma\\
&\xrightarrow{n}0.
\end{align*}
Further, combining the above two inequalities, by the arbitrariness of $\delta>0$, there holds
\begin{align}\label{R10}
\lim_{n\rightarrow\infty}\int_{\mathbb{R}^N}\big(I_\alpha\ast|u_{n,1}|^{2_*^\alpha}\big)|u_{n,2}|^{2_*^\alpha}dx=0.
\end{align}
Denote $a_2=\sqrt{a^2-a_1^2}$. Therefore, from \eqref{R8}, \eqref{R11}$-$\eqref{R10} and Lemma \ref{L1} we deduce
\begin{align}\label{R6}
m(a)&=\liminf_{n\to\infty}\mathcal{J}(u_n)\nonumber\\
&=\liminf_{n\to\infty}\bigg[\mathcal{J}(u_{n,1})+\mathcal{J}(u_{n,2})
+\frac12\int_{\mathbb{R}^N}|\nabla u_n|^2-|\nabla u_{n,1}|^2-|\nabla u_{n,2}|^2dx \nonumber\\
&\quad-\frac1{2\cdot2_*^\alpha}\int_{\mathbb{R}^N}\left(I_\alpha\ast|u_n|^{2_*^\alpha}\right)|u_n|^{2_*^\alpha}-
\left(I_\alpha\ast|u_{n,1}|^{2_*^\alpha}\right)|u_{n,1}|^{2_*^\alpha}
-\left(I_\alpha\ast|u_{n,2}|^{2_*^\alpha}\right)|u_{n,2}|^{2_*^\alpha}dx \nonumber\\
&\quad-\int_{\mathbb{R}^N}h(x)\left(u_n-u_{n,1}-u_{n,2}\right)dx \bigg]\nonumber\\
&=\liminf_{n\to\infty}\bigg[\mathcal{J}_\infty(u_{n,1})+\mathcal{J}(u_{n,2})
-\frac1{2_*^\alpha}\int_{\mathbb{R}^N}\left(I_\alpha\ast|u_{n,1}|^{2_*^\alpha}\right)|u_{n,2}|^{2_*^\alpha}dx\nonumber\\
&\quad +\frac12\int_{\mathbb{R}^N}|\nabla u_n|^2-|\nabla u_{n,1}|^2-|\nabla u_{n,2}|^2dx\bigg]\nonumber\\
&\geq\liminf_{n\to\infty}\mathcal{J}_\infty(u_{n,1})+\liminf_{n\to\infty}\mathcal{J}(u_{n,2})\nonumber\\
&\geq m_\infty(a_1)+m(a_2)\nonumber\\
&\geq m(a),
\end{align}
which implies $\liminf\limits_{n\to\infty}\mathcal{J}_\infty(u_{n,1})=m_\infty(a_1)$.
Set $\tilde{u}_{n,1}(\cdot)=u_{n,1}(\cdot+y_n)$, then $\liminf\limits_{n\to\infty}\mathcal{J}_\infty(\tilde{u}_{n,1})=m_\infty(a_1)$.
Since $\{\tilde{u}_{n,1}\}\subset H^1(\mathbb{B}_{2R_0}(0))$, the Rellich embedding theorem says $\tilde{u}_{n,1}\rightarrow\tilde{u}_1\in L^2(\mathbb{B}_{2R_0}(0))$. Clearly, $\tilde{u}_1\neq0$.
Further, it is easy to show $\mathcal{J}_\infty(\tilde{u}_1)=m_\infty(a_1)$, which is impossible due to Remark \ref{R1.2}.
If $\{y_n\}$ is bounded, assume $\sup\left\{|y_n|\right\}\leq\tau$, then $u_{n,2}(x)\rightarrow0$ a.e.~in $\mathbb{R}^N$.
Similar to the proof of \eqref{R6}, we have
\begin{align}\label{R7}
m(a)&=\liminf_{n\to\infty}\mathcal{J}(u_n)\nonumber\\
&\geq\liminf_{n\to\infty}\mathcal{J}(u_{n,1})+\liminf_{n\to\infty}\mathcal{J}_\infty(u_{n,2})\nonumber\\
&\geq m(a_1)+m_\infty(a_2)\nonumber\\
&\geq m(a),
\end{align}
which implies $\liminf\limits_{n\to\infty}\mathcal{J}(u_{n,1})=m(a_1)$.
Since $\{u_{n,1}\}\subset H^1(\mathbb{B}_{2R_0+\tau}(0))$, the Rellich embedding theorem says $u_{n,1}\rightarrow \hat{u}_1\in L^2(\mathbb{B}_{2R_0+\tau}(0))\backslash\{0\}$. Then we derive $\mathcal{J}(\hat{u}_1)=m(a_1)$.
Take $\eta\in C^\infty(\mathbb{R}^N,[0,1])$ such that $\eta(x)=0$ if $|x|\leq2R_0+\tau$ and $\eta(x)=1$ if $|x|\geq3R_0+\tau$.
Set $u_\varepsilon(\cdot)=\frac{a_2\varepsilon^{\frac N2}u_0(\varepsilon\cdot)\eta(\cdot)}{|\eta_\varepsilon u_0|_2}$ with $u_0$ given in the proof of Lemma \ref{L1} and $\eta_\varepsilon(\cdot)=\eta(\frac\cdot\varepsilon)$, by direct computations, there hold for $\varepsilon>0$ small enough
\begin{align}\label{R5}
\begin{cases}
\vspace{0.10cm}
\displaystyle |u_0|_2^2-O(\varepsilon^N)=\left|\eta_\varepsilon u_0\right|_2^2\leq |u_0|_2^2, \\
\vspace{0.10cm}
\displaystyle \int_{\mathbb{R}^N}|\nabla u_\varepsilon|^2dx=\frac{|u_0|_2^2}{|\eta_\varepsilon u_0|_2^2}\int_{\mathbb{R}^N}|\nabla u_{a_2}^\varepsilon|^2dx+O(\varepsilon^N), \\
\vspace{0.10cm}
\displaystyle \int_{\mathbb{R}^N}h(x)u_\varepsilon dx
=\frac{|u_0|_2}{|\eta_\varepsilon u_0|_2}\int_{\mathbb{R}^N}h(x)u_{a_2}^\varepsilon dx-O(\varepsilon^{\frac N2}), \\
\vspace{0.10cm}
\displaystyle \int_{\mathbb{R}^N}\big(I_\alpha\ast|u_\varepsilon|^{2_*^\alpha}\big)|u_\varepsilon|^{2_*^\alpha}dx
\geq\left(\frac{|u_0|_2^2}{|\eta_\varepsilon u_0|_2^2}\right)^{2_*^\alpha}
\int_{\mathbb{R}^N}\big(I_\alpha\ast|u_{a_2}^\varepsilon|^{2_*^\alpha}\big)|u_{a_2}^\varepsilon|^{2_*^\alpha}dx-O(\varepsilon^{\frac{N+\alpha}2}),
\end{cases}
\end{align}
where $u_{a_2}^\varepsilon(\cdot)=\frac{a_2\varepsilon^{\frac N2}u_0(\varepsilon\cdot)}{|u_0|_2}$.
It is clear that $|\hat{u}_1|_2=a_1$ and $|u_\varepsilon|_2=a_2$.
Then, observe $\mbox{spt}~\hat{u}_1\cap \mbox{spt}~u_\varepsilon=\emptyset$ and $|\hat{u}_1+u_\varepsilon|_2=a$,
by \eqref{R4}, \eqref{R5} and ($\mathcal{H}$) we conclude, for $\varepsilon>0$ small~enough,
\begin{align*}
m(a)
&\leq\mathcal{J}(\hat{u}_1+u_\varepsilon)\leq\mathcal{J}(\hat{u}_1)+\mathcal{J}(u_\varepsilon) \\
&\leq m(a_1)+\frac{a_2^2\varepsilon^2}{|u_0|_2^2}\int_{\mathbb{R}^N}|\nabla u_0|^2dx
-\frac1{2\cdot2_*^\alpha}\left(\frac {a_2}{|u_0|_2}\right)^{2\cdot2_*^\alpha}
-\frac{a_2\varepsilon^{\frac N2}}{|u_0|_2}\int_{\mathbb{B}_{\frac {\bar{r}}\varepsilon}^c(0)}
\frac{u_0(\varepsilon x)}{|x|^\beta}dx+O(\varepsilon^{\frac{N}2}) \\
&\leq m(a_1)-\frac1{2\cdot2_*^\alpha}\left(\frac {a_2^2}{S_\alpha}\right)^{2_*^\alpha}
+\frac{a_2^2\varepsilon^2}{|u_0|_2^2}\int_{\mathbb{R}^N}|\nabla u_0|^2dx
-\frac{a_2\varepsilon^{\beta-\frac N2}}{|u_0|_2}\int_{\mathbb{B}_{\bar{r}}^c(0)}\frac{u_0(x)}{|x|^\beta}dx
+O(\varepsilon^{\frac{N}2}) \\
&<m(a_1)+m_\infty(a_2),
\end{align*}
which contradicts \eqref{R7}. Hence, compactness occurs.
Since $\{u_n\}$ is bounded, there exists $u\in H^1(\mathbb{R}^N)$ such that, up to a subsequence,
\begin{align}\label{R2}
u_n\rightharpoonup u\ \mbox{in}\ H^1(\mathbb{R}^N), \ \ \ \ \
u_n\rightarrow u\  \mbox{in}\ L^2_{loc}(\mathbb{R}^N)\ \ \ \ \mbox{and}\ \ \ \
u_n(x)\rightarrow u(x)\ a.e.\ \mbox{in}\ \mathbb{R}^N.
\end{align}

If $\{z_n\}$ is unbounded. Suppose $|z_n|\to\infty$ in the sense of subsequence.
Observe that $\{u_n(\cdot+z_n)\}$ is bounded, then $u_n(\cdot+z_n)\rightharpoonup u_0$ in $H^1(\mathbb{R}^N)$ and
$u_n(\cdot+z_n)\rightarrow u_0$ in $L^2_{loc}(\mathbb{R}^N)$ up to a subsequence.~As a consequence, by \eqref{R1}, $|u_0|_2=a$ and then $u_n(\cdot+z_n)\rightarrow u_0$ in $L^2(\mathbb{R}^N)$.
It is easy to show, as $n\rightarrow\infty$,
\begin{align*}
\int_{\mathbb{R}^N}u_0(x-z_n)\varphi(x)dx=\int_{\mathbb{R}^N}u_0(x)\varphi(x+z_n)dx\xrightarrow{n}0,\ \ \
\forall\ \varphi\in C_0^\infty(\mathbb{R}^N),
\end{align*}
which and the density of $C_0^\infty(\mathbb{R}^N)$ in $H^1(\mathbb{R}^N)$ imply $u_0(x-z_n)\rightharpoonup0$ in $L^2(\mathbb{R}^N)$.
Then, by $h\in L^2(\mathbb{R}^N)$,
\begin{align*}
\lim_{n\rightarrow\infty}\int_{\mathbb{R}^N}h(x+z_n)u_n(x+z_n)dx=\lim_{n\rightarrow\infty}\int_{\mathbb{R}^N}h(x)u_0(x-z_n)dx=0.
\end{align*}
Moreover, it is obvious that $\|u_n(\cdot+z_n)\|_*\rightarrow\|u_0\|_*$.
Combining the above arguments, by the weakly lower semicontinuity of the norm and that $m_\infty(a)$ has no minimizer, we conclude
\begin{align*}
m(a)&=\lim_{n\to\infty}\mathcal{J}(u_n) \\
&=\lim_{n\to\infty}\bigg[\frac12\int_{\mathbb{R}^N}|\nabla u_n(x+z_n)|^2dx
-\frac1{2\cdot2_*^\alpha}\int_{\mathbb{R}^N}\left(I_\alpha\ast|u_n(\cdot+z_n)|^{2_*^\alpha}\right)|u_n(x+z_n)|^{2_*^\alpha}dx\\
&\quad-\int_{\mathbb{R}^N}h(x+z_n)u_n(x+z_n)dx\bigg] \\
&\geq\mathcal{J}_\infty(u_0)\\
&>m_\infty(a),
\end{align*}
which contradicts with Lemma \ref{L1}-(1). Consequently, $\{z_n\}$ is bounded.
Further, it follows from \eqref{R1} and \eqref{R2} that $|u|_2=a$ and then $u_n\rightarrow u$ in $L^2(\mathbb{R}^N)$.
In view of this, it is easy to verify
\begin{align*}
m(a)=\lim_{n\to\infty}\mathcal{J}(u_n)\geq\mathcal{J}(u)\geq m(a),
\end{align*}
which implies $\mathcal{J}(u)=m(a)$.
It is clear that $\mathcal{J}(u)\geq\mathcal{J}(|u|)\geq m(a)$ since $h(\cdot)\geq0$, i.e. $\mathcal{J}(|u|)=m(a)$.
Hence, we can assume $u(\cdot)\geq0$. By the Lagrange multiplier rule, there exists some $\lambda_a\in \mathbb{R}$ such that
\begin{align*}
-\Delta u+\lambda_a u=\left(I_\alpha\ast|u|^{2_*^\alpha}\right)|u|^{2_*^\alpha-2}u+h(x)\ \ \ \mbox{in}\ \mathbb{R}^N.
\end{align*}
Thanks to the strong maximum principle, we conclude $u(\cdot)>0$.
Based on the above arguments, $u$ is a positive ground state normalized solution to Eq.~\eqref{CH3}, with the Lagrange multiplier
\begin{align*}
\lambda_a&=\frac1{a^2}\left[-\int_{\mathbb{R}^N}|\nabla u|^2dx+\int_{\mathbb{R}^N}\left(I_\alpha\ast|u|^{2_*^\alpha}\right)|u|^{2_*^\alpha}dx+\int_{\mathbb{R}^N}h(x)udx\right] \\
&=\frac1{a^2}\left[-2\mathcal{J}(u)+\frac\alpha{N+\alpha}\int_{\mathbb{R}^N}\left(I_\alpha\ast|u|^{2_*^\alpha}\right)|u|^{2_*^\alpha}dx
-\int_{\mathbb{R}^N}h(x)udx\right] \\
&>\frac{a^{2(2_*^\alpha-1)}}{2_*^\alpha S_\alpha^{2_*^\alpha}}-\frac{|h|_2}{a},
\end{align*}
where we used the H$\ddot{\mbox{o}}$lder inequality and $h\in L^2(\mathbb{R}^N)$. Thus we complete the proof of Theorem \ref{T1}.
\end{prf1}\\

Next, we turn to prove the nonexistence of nonnegative minimizer for $m(a)$ in case $h(\cdot)\leq0$. \\

\begin{prf2}
Similar to the proof of \eqref{R4}, noting $-u_a^\varepsilon(\cdot)\in\mathcal{S}(a)$ for any $\varepsilon>0$, we deduce
\begin{align*}
m(a)&\leq\lim_{\varepsilon\rightarrow0^+}\mathcal{J}(-u_a^\varepsilon)\\
&=\lim_{\varepsilon\rightarrow0^+}\bigg[\frac{a^2\varepsilon^2}{2|u_0|_2^2}\int_{\mathbb{R}^N}|\nabla u_0|^2dx
-\frac1{2\cdot2_*^\alpha}\left(\frac a{|u_0|_2}\right)^{2\cdot2_*^\alpha}
+\frac{a\varepsilon^{\frac N2}}{|u_0|_2}\int_{\mathbb{R}^N}h(x)u_0(\varepsilon\cdot) dx\bigg] \\
&\leq\lim_{\varepsilon\rightarrow0^+}\bigg[ \frac{a^2\varepsilon^2}{2|u_0|_2^2}\int_{\mathbb{R}^N}|\nabla u_0|^2dx
-\frac1{2\cdot2_*^\alpha}\left(\frac {a^2}{S_\alpha}\right)^{2_*^\alpha}\bigg]\\
&=-\frac1{2\cdot2_*^\alpha}\left(\frac {a^2}{S_\alpha}\right)^{2_*^\alpha}.
\end{align*}
Then, if there is $u\in \mathcal{S}(a)$ such that $u(\cdot)\geq0$ and $\mathcal{J}(u)=m(a)$, it follows from the definition of~$S_\alpha$~that
\begin{align*}
0&<\frac12\int_{\mathbb{R}^N}|\nabla u|^2dx-\int_{\mathbb{R}^N}h(x)u dx \\
&\leq\frac1{2\cdot2_*^\alpha}\int_{\mathbb{R}^N}\left(I_\alpha\ast|u|^{2_*^\alpha}\right)|u|^{2_*^\alpha}dx
-\frac1{2\cdot2_*^\alpha}\left(\frac {a^2}{S_\alpha}\right)^{2_*^\alpha} \\
&\leq0,
\end{align*}
which is a contradiction. Therefore, $m(a)$ has no nonnegative minimizer in case $h$ is non-positive.
\end{prf2}



\end{document}